\pdfoutput=1
\documentclass[12pt]{article}
\usepackage[english]{babel}
\usepackage{authblk}
\usepackage{numprint}
\usepackage{a4}
\usepackage{cite}
\usepackage{graphicx}
\usepackage{rotating}
\usepackage{hyperref}
\usepackage{xspace}
\usepackage{subfigure} 
\usepackage[T1]{fontenc}
\usepackage[latin1]{inputenc}
\def\3{\ss}

\def\co#1{\texttt{#1}\xspace}
\title{How to get meaningful and correct results from your finite element model}
\author{Martin B\"aker}
\affil{Institut f\"ur Werkstoffe, Technische Universit\"at
  Braunschweig, Langer Kamp 8, D-38106 Braunschweig, martin.baeker@tu-bs.de}
\begin{document}
\maketitle

\begin{abstract}
This document gives guidelines to set up, run, and postprocess correct
simulations with the finite element method. It is not an introduction
to the method itself, but rather a list of things to check and
possible mistakes to watch out for when doing a finite element simulation.
\end{abstract}

The finite element method (FEM) is probably the most-used simulation
technique in engineering.  Modern finite-element software makes doing
FE simulations easy~-- too easy, perhaps. Since you have a nice
graphical user interface that guides you through the process of
creating, solving, and postprocessing a finite element model, it may
seem as if there is no need to know much about the inner workings of a
finite element program or the underlying theory. However, creating a
model without understanding finite elements is similar to flying an
airplane without a pilot's license. You may even land somewhere
without crashing, but probably not where you intended to.

This document is not a finite element introduction, see, for example,
\cite{Dhondt,baker2002numerische,bathe2006} for that.  It is a
guideline to give you some ideas how to correctly set up, solve and
postprocess a finite element model. The techniques described here were
developed working with the program Abaqus \cite{abaqus}; however, most
of them should be easily transferable to other codes. I have not
explained the theoretical basis for most of them; if you do not
understand why a particular consideration is important, I recommend to
study finite element theory to find out.

\section{Setting up the model}
\renewcommand*{\theenumi}{\thesubsection-\arabic{enumi}}
\renewcommand*{\theenumii}{\theenumi.\arabic{enumii}}

\subsection{General considerations}
These considerations are not restricted to finite element models, but
are useful for any complex simulation method.

\begin{enumerate}
\item Even if you just need some number for your design~-- the main
  goal of an FEA is to \emph{understand} the system. Always
  design your simulations so that you can at least qualitatively
  understand the results. \emph{Never} believe the result of a
  simulation without thinking about its plausibility. 
\item Define the goal of the simulation as precisely as
  possible. Which question is to be answered? Which quantities are to
  be calculated? Which conclusions are you going to draw from the
  simulation? Probably the most common error made in FE
  simulations is setting up a simulation without having a clear goal
  in mind. Be as specific as possible. Never set up a model ``to
  see what happens'' or ``to see how stresses are
  distributed''. 
\item Formulate your expectations for the simulation result beforehand
  and make an educated guess of what the results should be. If
  possible, estimate at least some quantities of your simulation using
  simplified assumptions. This will make it easier to spot problems
  later on and to improve your understanding of the system you are
  studying. 
\item Based on the answer to the previous items, consider which effects
  you actually have to simulate. Keep the model as simple as
  possible. For example, if you only need to know whether a yield
  stress is exceeded somewhere in a metallic component, it is much easier to
  perform an elastic calculation and check the von Mises stress in the
  postprocessor (be wary of extrapolations, see \ref{item:extrapolate}) than to include plasticity in the model.
\item What is the required precision of your calculation? Do you need
  an estimate or a precise number? (See also~\ref{item:inputPrecision} below.) 
\item If your model is complex, create it in several steps. Start with
  simple materials, assume frictionless behaviour
  etc. Add complications step by step. Setting up the model in
  steps has two advantages: (i) if errors 
  occur, it is much easier to find out what caused them; (ii)
  understanding the behaviour of the system is easier this way because
  you understand which addition caused which change in the model behaviour.
  Note, however, that checks you made in an early stage (for example
  on the mesh density) may have to be repeated later.
\item Be careful with units. Many FEM programs (like ABAQUS) are
  inherently unit-free~-- they assume that all numbers you give can be
  converted without additional conversion factors. You cannot define
  you model geometry in millimeter, but use SI units without prefixes everywhere
  else. Be especially careful in thermomechanical simulations due to
  the large number of different physical quantities needed there. And
  of course, be also careful if you use antiquated units like inch,
  slug, or BTU.
\end{enumerate}

\subsection{Basic model definition}
\begin{enumerate}
\item Choose the correct type of simulation (static, quasi-static,
  dynamic, coupled etc.). Dynamic simulations
require the presence of inertial forces (elastic waves, changes in
kinetic energies). If inertial
forces are irrelevant, you should use static simulations.
\item As a rule of thumb, a simulation is static or quasi-static if
  the excitation frequency is less than $1/5$ of the lowest natural
  frequency of the structure \cite{caeai}.
\item In a dynamic analysis, damping may be required to avoid
  unrealistic multiple
  reflections of elastic waves that may affect the results \cite{caeai}.
\item \label{it:explicitScaling}
Explicit methods are inherently dynamic. In some cases, explicit
  methods may be used successfully for quasi-static problems to avoid
  convergence problems (see~\ref{it:explicitConvergence} below). 
If you use mass scaling in your explicit quasi-static analysis,
carefully check that the scaling parameter does not affect your
solution. Vary the scaling factor (the nominal density) to ensure that
the kinetic energy in the model remains small \cite{prior1994}.  
\item In a static or quasi-static analysis, make sure that all parts
  of the model are constrained so that no rigid-body movement is
  possible. (In a contact problem, special stabilization techniques
  may be available to ensure correct behaviour before contact is established.)
\item If you are studying a coupled problem (for example thermo-mechanical) think about the correct form of
  coupling. If stresses and strains are affected by temperature but
  not the other way round, it may be more efficient to first calculate
  the thermal problem and then use the result to calculate thermal
  stresses. A full coupling of the thermal and mechanical problem is
  only needed if temperature affects stresses/strains (e.\,g., due to thermal
  expansion or temperature-dependent material problems) and if
  stresses and strains also affect the thermal problem (e.\,g., due to plastic
  heat generation or the change in shape affecting heat conduction).
\item\label{item:timeSteps} Every FE program uses discrete time steps
  (except for a static, linear analysis, where no time incrementation
  is needed). This may affect the
  simulation.  If, for example, the temperature changes during a time
  increment, the material behaviour may strongly differ between the
  beginning and the end of the increment (this often occurs in creep
  problems where the properties change drastically with
  temperature). Try different maximal time increments and make sure
  that time increments are sufficiently small so that these effects
  are small.
\item\label{item:nlgeom} Critically check whether non-linear geometry is required. As a
  rule of thumb, this
  is almost always the case if strains exceed 5\%. If loads are
  rotating with the structure (think of a 
  fishing rod that is loaded in bending initially, but in tension
  after it has started to deform),
  the geometry is usually non-linear. If in doubt, critically compare a
  geometrically linear
  and non-linear simulation.
\end{enumerate}

\subsection{Symmetries, boundary conditions and loads}\label{sec:symmetries}
\begin{enumerate}
\item Exploit symmetries of the model. In a plane 2D-model, think about
  whether plane stress, plane strain or generalized plane strain is
  the appropriate symmetry. (If thermal stresses are relevant, plane
  strain is almost always wrong because thermal expansion in the
  3-direction is suppressed, causing large thermal stresses. Note that
  these 33-stresses may affect other stress components as well, for example, due
  to von Mises plasticity.) Keep in mind that the loads and the
  deformations must conform to the same symmetry.

\item Check boundary conditions and constraints. After calculating the
  model, take the time to ensure that nodes were
  constrained in the desired way  in the postprocessor. 
\item Point loads at single nodes may cause unrealistic stresses in
  the adjacent elements. Be especially careful if the material or the
  geometry is non-linear. If in doubt, distribute the load over
  several elements (using a local mesh refinement if necessary).
\item If loads are changing direction during the calculation,
  non-linear geometry is usually required, see \ref{item:nlgeom}. 
\item\label{item:timestep2} The discrete time-stepping of the solution
  process may also be important in loading a structure. If, for
  example, you abruptly change the heat flux at a certain point in
  time, discrete time stepping may not capture the exact point at
  which the change occurs, see fig.~\ref{fig:timeStep}. (Your software may use some
  averaging procedure to alleviate this.) Define load
  steps or use other methods to ensure that the time of the abrupt
  change actually corresponds to a time step in the simulation. This
  may also improve convergence because it allows to control the
  increments at the moment of abrupt change, see also~\ref{item:rapidChange}

\begin{figure}
\includegraphics[width=0.7\textwidth]{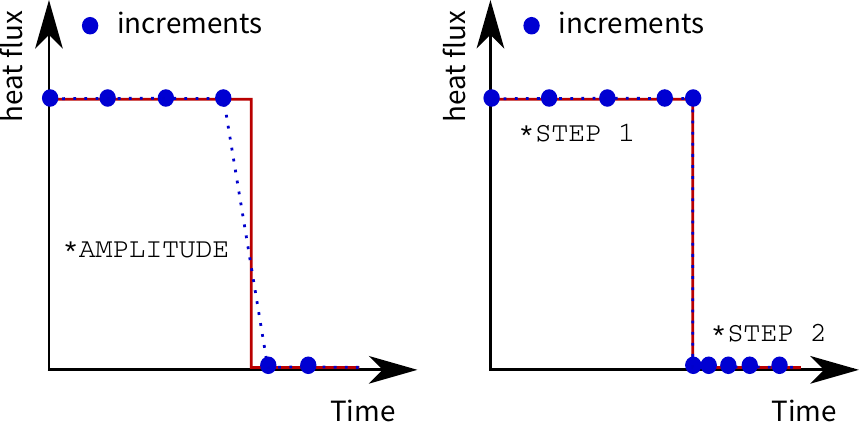}
\caption{Discrete time steps may affect the load of a structure. If,
  for example, a heat flux abruptly changes, the change may actually
  only become relevant at a later moment (left). Use a new load step (or
  another method) to ensure that you capture the correct moment (right).}
\label{fig:timeStep}
\end{figure}
\end{enumerate}

\subsection{Input data}\label{sec:eingangsgroessen}
\begin{enumerate}
\item\label{item:inputPrecision}
A simulation cannot be more precise than its input data allow. This is
especially true for the material behaviour. Critically consider how
precise your material data really are. How large are the
uncertainties? If in doubt, vary material parameters to see how
results are affected by the uncertainties.
\item Be careful when combining material data from different sources
  and make sure that they are referring to identical materials. In
  metals, don't forget to check the influence of heat treatment; in
  ceramics, powder size or the processing route may affect the
  properties; in polymers, the chain length or the content of
  plasticizers is important \cite{rosler2007}. Carefully document your sources for
  material data and check for inconsistencies.
\item Be careful when extrapolating material data. If data have been
  described using simple relations (for example a Ramberg-Osgood law
  for plasticity), the real behaviour may strongly deviate from this.
\item Keep in mind that your finite element software usually cannot
  extrapolate material data beyond the values given. If plastic
  strains exceed the maximum value specified, usually no further
  hardening of the material will be considered. The same holds, for
  example,  for
  thermal expansion coefficients which usually increase with
  temperature. Using different ranges in different materials may thus
  cause spurious thermal stresses. Fig.~\ref{fig:CTE} shows an example
\begin{figure}
\includegraphics[width=0.6\textwidth]{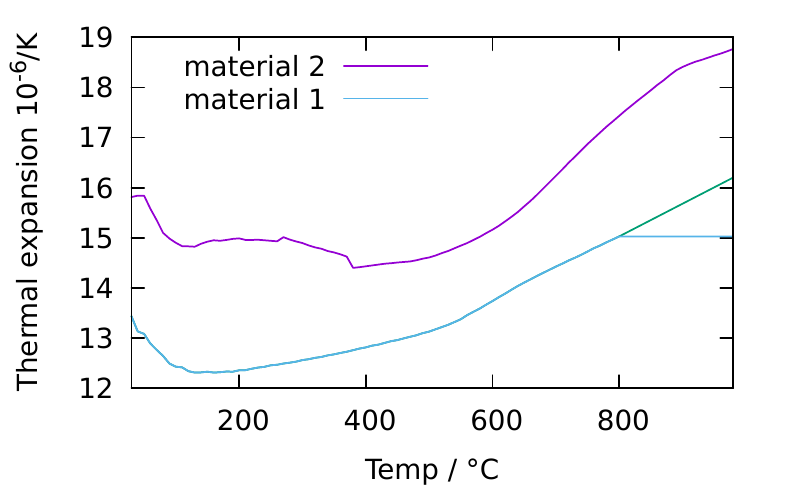}
\caption{Coefficient of thermal expansion for two
  materials. Experimental data for material 1 end at a temperature of
  \numprint[{}^\circ C]{800}, data for material 2 extend to a higher temperature. Most finite element programs assume a
  constant value beyond the final data point (blue line). If data are
  not extrapolated, the mismatch of the expansion coefficient  between
  the materials is overestimated, causing large thermal stresses.}
\label{fig:CTE}
\end{figure}
\item If material data are given as equations, be aware that
  parameters may not be unique. Frequently, data can be fitted using
  different parameters. As an illustration, plot the simple hardening law
 $A+B\epsilon^n$ with values $(130,
100,0.5)$ and $(100,130,0.3)$ for $(A,B,n)$, see fig.~\ref{fig:nonUnique}.
\begin{figure}
\includegraphics[width=0.4\textwidth]{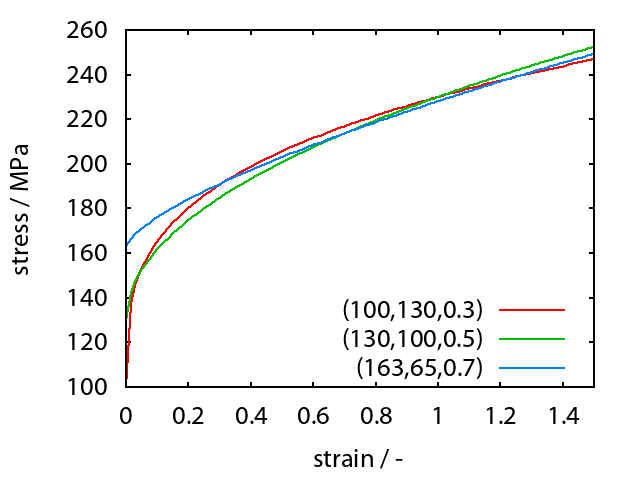}
\caption{Different values of the parameters $A$,$B$, and
  $n$ in a flow stress law  $\sigma=A+B\epsilon^n$  result in very similar curves.}
\label{fig:nonUnique}
\end{figure}
Your simulation results may be indifferent to some changes in the
parameters because of this.
\item If it is not possible to determine material behaviour precisely,
  finite element simulations may still help to understand how
  the material behaviour affects the system. Vary parameters in
  plausible regions and study the answer of the system.
\item Also check the precision of external loads. If loads are not known precisely, use
  a conservative estimate.
\item
 Thermal loads may be especially
  problematic because heat transfer coefficients or surface
  temperatures may be difficult to measure. Use the same
  considerations as for materials. 
\item If you vary parameters (for example the geometry of your
  component or the material), make sure that you correctly consider
  how external loads are changed by this. If, for example, you specify
  an external load as a pressure, increasing the surface also
  increases the load. If you change the thermal conductivity of your
  material, the total heat flux through the structure will change; you
  may have to specifiy the thermal load accordingly.
\item Frictional behaviour and friction coefficients are also
  frequently unknown. Critically check the parameters you use and also
  check whether the friction law you are using is correct~-- not all
  friction is Coulombian.
 \item If a small number of parameters are unknown, you can try to vary them
   until your simulation matches experimental data, possibly using a
   numerical optimization method. (This is the
   so-called inverse parameter identification \cite{baker2013}.) Be aware that the
   experimental data used this way cannot be used to validate your model
   (see section~\ref{sec:verifikation}).
\end{enumerate}
\subsection{Choice of the element type}
\emph{Warning:} Choosing the element type is often the crucial step in
creating a finite element model. Never accept the default choice of
your program without thinking about it.\footnote{The only acceptable
  exception may be a simple linear-elastic simulation if your program
  uses second-order elements. But if all you do is linear elasticity,
  this article is probably not for you.} Carefully check which types are
available and make sure you understand how a finite element simulation
is affected by the choice of element type. You should understand the
concepts of element order and integration points (also known as
Gau\ss\ points) and know the most common errors caused by an
incorrectly chosen element type (shear locking, volumetric locking,
hourglassing \cite{armero2000,baker2002numerische}).
 
The following points give some guidelines for the correct choice :
\begin{enumerate}
\item If your problem is linear-elastic, use second-order
  elements. Reduced integration may save computing time without
  strongly affecting the results.
\item Do not use fully-integrated first order elements if bending
  occurs in your structure (shear locking). Incompatible mode elements
  may circumvent this problem, but their performance strongly depends
  on the element shape \cite{bathe2006}.
\item If you use first-order elements with reduced integration, check
  for hourglassing. Keep in mind that hourglassing may occur only in
  the interior of a three-dimensional structure where seeing it is not
  easy. Exaggerating the displacements may help in visualizing
  hourglassing. Most programs use numerical techniques to suppress
  hourglass modes; however, these may also
  affect results due to artificial damping. Therefore, also check the
  energy dissipated by this artificial damping and make sure that it
  is small compared to other energies in the model.
\item\label{item:contactSecond} In contact problems, first-order elements may improve
  convergence because if one corner and one edge node are in contact,
  the second-order interpolation of the element edge causes overlaps,
  see fig.~\ref{fig:nodeRelease}. This may especially cause problems
  in a crack-propagation simulation with a node-release scheme
  \cite{Krueger2004,Baeker2018}.
\begin{figure}
\includegraphics[width=0.3\textwidth]{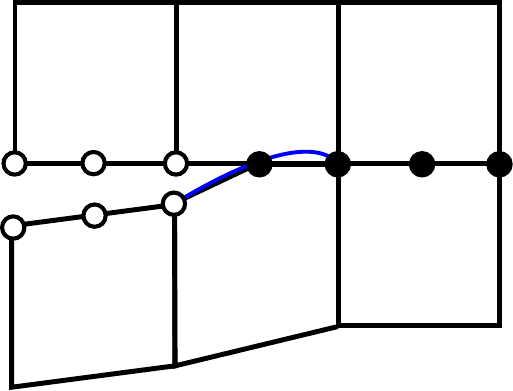}
\caption{Contact between second-order elements. If an edge and one
  corner node of an element are in contact with a surface, the
  quadratic interpolation of the edge shape causes a penetration of
  the elements. Adapted from \cite{Baeker2018}.}
\label{fig:nodeRelease}
\end{figure}
\item Discontinuities in stresses or strains may be captured better
  with first-order elements in some circumstances.
\item If elements distort strongly, first-order elements may be better
  than second-order elements.
\item\label{item:triang}
 Avoid triangular or tetrahedral first-order elements since they
  are much too stiff, especially in bending. If you have to use these
  elements (which may be necessary in a large model with complex
  geometry), use a very fine mesh and carefully check for mesh
  convergence. Think about whether partitioning your model and meshing
  with quadrilateral/hexahedral elements (at least in critical
  regions) may be worth the effort.  Fig.~\ref{fig:foam} shows an
  example where a very complex geometry has to be meshed with
  tetrahedral elements. Although the mesh looks reasonably fine, the system answer
  with linear elements is much too stiff.

\begin{figure}
\includegraphics[width=0.6\textwidth]{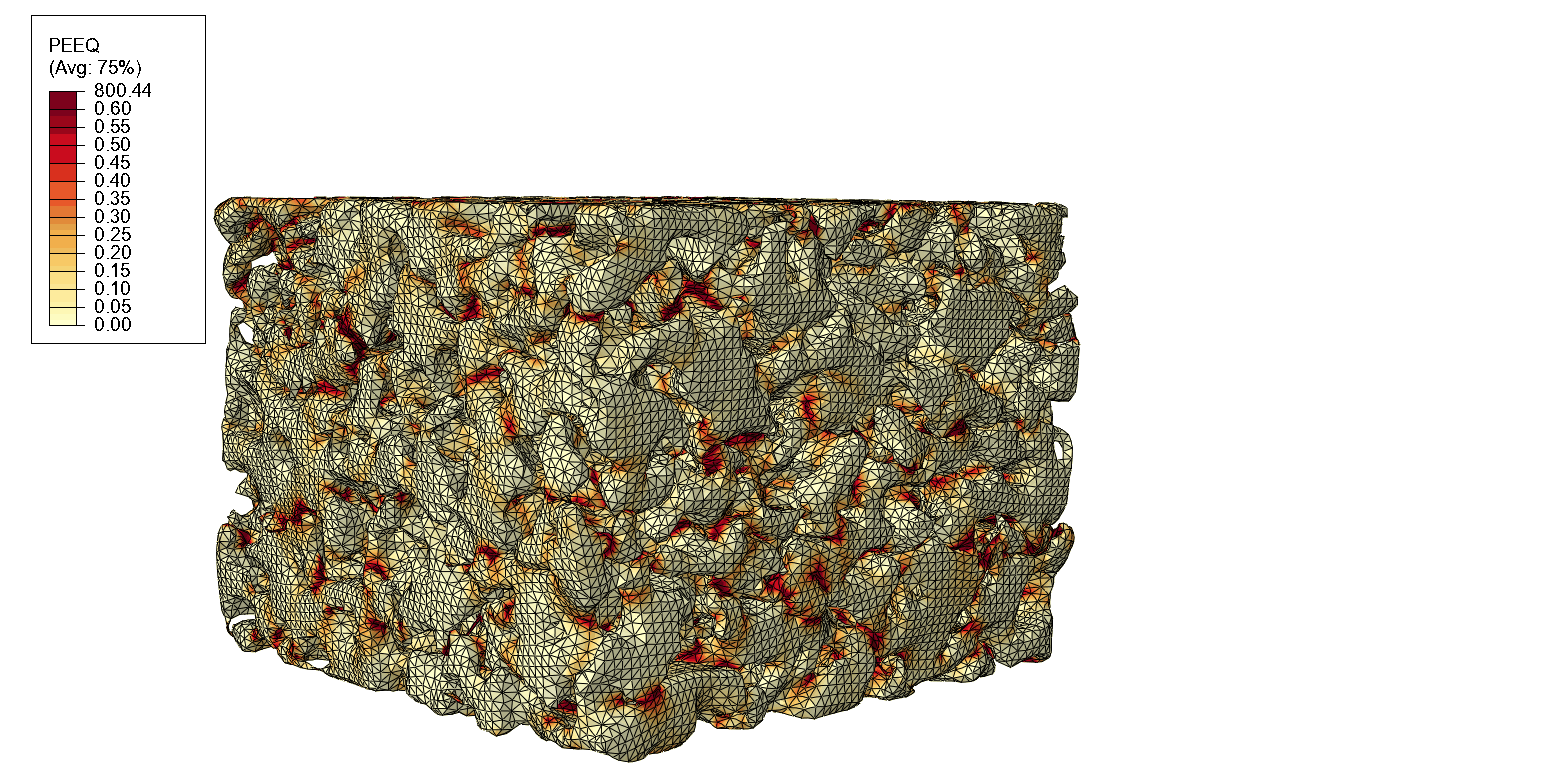}
\includegraphics[width=0.38\textwidth]{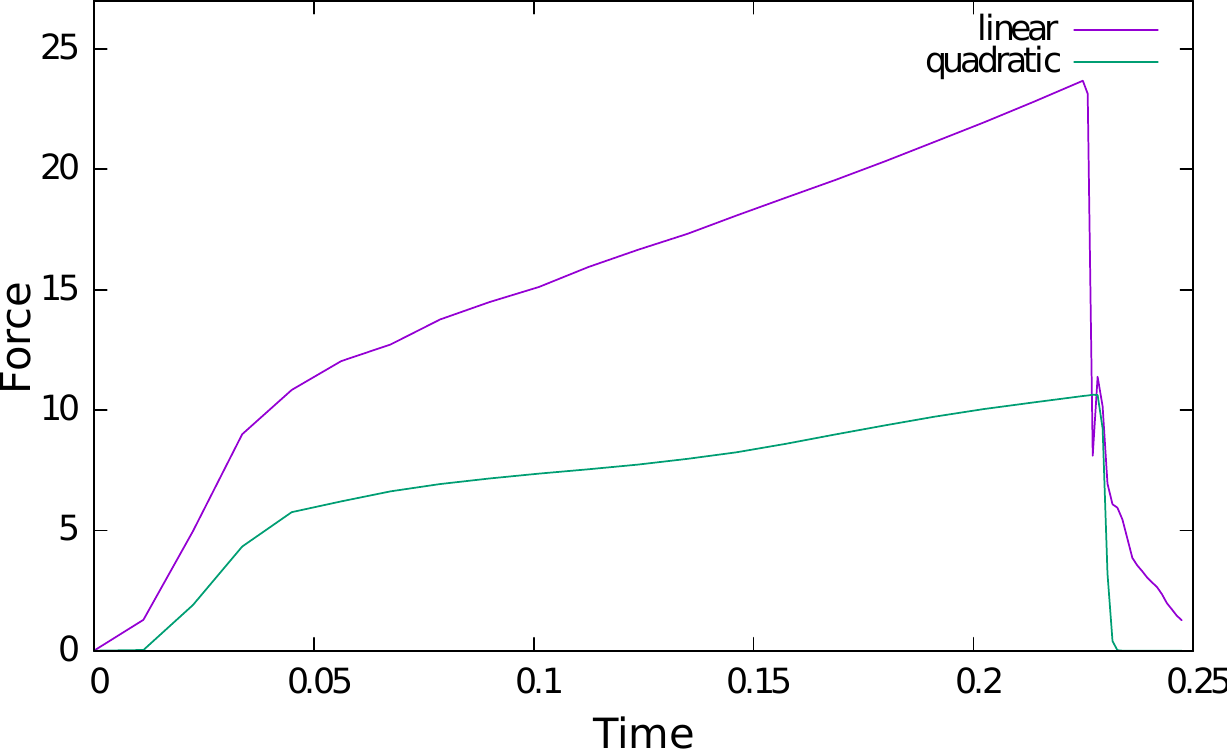} 
\caption{Simulation of the compression of a metallic foam with a
  tetrahedral mesh. As the force-time curve shows, the result
  strongly differs between linear and quadratic elements although the
  mesh looks rather fine and comprises more than 700000 elements. Note
that the simulation is displacement-controlled so evaluating forces
reveals whether the model is too stiff, see~\ref{item:meshconverge}.}
\label{fig:foam}
\end{figure}
\item If material behaviour is incompressible or almost
  incompressible, use hybrid elements to avoid volumetric locking. They may also be useful if
  plastic deformation is large because (metal) plasticity is also volume
  conserving. 
\item Do not mix elements with different order. This can cause
  overlaps or gaps forming at the interface (possibly not shown by
  your postprocessor) even if there are no hanging nodes (see
  fig.~\ref{fig:mixing}). If you have to use 
  different order of elements in different regions of your model, tie
  the interface between the regions using a surface constraint. Be
  aware that this interface may cause a discontinuity in the stresses
  and strains due to different stiffness of the element types.
\begin{figure}
\includegraphics[width=0.4\textwidth]{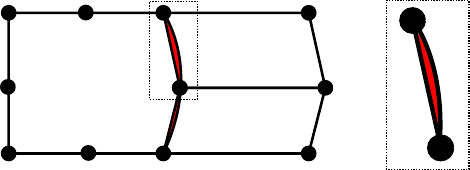}
\caption{Mixing of elements with different order (left: one
  second-order element, right: two first-order elements) is not
  allowed as it would lead to unphysical deformation (gaps or
  overlaps): The quadratic element uses a second-order function to
  calculate the elements edge, the two linear elements assume a
  piecewise linear shape.}
\label{fig:mixing}
\end{figure}
\item In principle, it is permissible to mix reduced and fully
  integrated elements of the same order. However, since they differ in
  stiffness, spurious stress or strain discontinuities may result.
\item If you use shell or beam elements or similar, make sure to use
  the correct formulation. Shells and membranes look similar, but
  behave differently. Make sure that you use
  the correct mathematical formulation; there are a large number of
  different types of shell or beam elements with different behaviour. 
\end{enumerate}

\subsection{Generating a mesh}\label{sec:mesh}

\begin{enumerate}
\item If possible, use quadrilateral/hexahedral elements. Meshing
  3D-structures this way may be laborious, but it is often worth the
  effort (see also~\ref{item:triang}).
\item A fine mesh is needed where gradients in stress and strain are large.
\item A preliminary simulation with a coarse mesh may help to identify
  the regions where a greater mesh density is required.
\item Keep in mind that the required mesh density depends on the
  quantities you want to extract and on the required precision. For example, displacements are often
  calculated more precisely than strains (or stresses) because strains involve
  derivatives, i.e.\ the
  differences in displacements between nodes. 
\item \label{item:meshconverge}
A mesh convergence study can be used to check whether the model
  behaves too stiff (as is often the case for fully integrated
  first-order elements, see fig.~\ref{fig:foam}) or too soft (which happens with
  reduced-integration elements). Be careful in evaluating this
  study: If your model is load-controlled, evaluate displacements or
  strains to check for convergence, if it is strain-controlled,
  evaluate forces or stresses. (Stiffness relates forces to
  displacements, so to check for stiffness you need to check both.) If
  you use, for example, displacement  
  control, displacements are not sensitive to the actual stiffness
  of your model since you prescribe the displacement.
\item Check shape and size of the elements. Inner angles should not
  deviate too much from those of a regularly shaped element. Use the
  tools provided by your software to highlight critical elements. Keep
  in mind that critical regions may be situated inside a 3D-component
  and may not be directly visible. Avoid badly-shaped elements
  especially in region where high gradients occur and in regions of interest.
\item If you use local mesh refinement, the transition between regions
  of different element sizes should be smooth. As a rule of thumb,
  adjacent elements should not differ by more than a factor of 2--3 in
  their area (or volume). If the transition is too abrupt, spurious stresses may
  occur in this region because a region that is meshed finer is
  usually less stiff. Furthermore,  the fine mesh may be constrained by the
  coarser mesh. (As an extreme case, consider a finely meshed
quadratic region that is bounded by only four first-order elements~-- in this
case, the region as a whole can only deform as a parallelogram, no
matter how fine the interior mesh is.)
\item Be aware that local mesh refinement may strongly affect the
  simulation time in an explicit simulation because the stable time
  increment is determined by the size of the smallest element in the
  structure. A single small or badly shaped element can drastically
  increase the simulation time.
\item If elements are distorting strongly, remeshing may improve the
  shape of the elements and the solution quality. For this, solution
  variables have to be interpolated from the old to the new mesh. This
  interpolation may dampen strong gradients or local extrema. Make
  sure that this effect is sufficiently small by comparing the
  solution before and after the remeshing in a contour plot and at the
  integration points.
\item \label{it:distort}
Another way of dealing with strong mesh distortions is to start
  with a mesh that is initially distorted and becomes more regular
  during deformation. This method usually requires some
  experimentation, but it may yield good solutions without the
  additional effort of remeshing.
\end{enumerate}

\subsection{Defining contact problems}

\begin{enumerate}
\item Correctly choose master and slave surfaces in a master-slave
  algorithm. In general, the stiffer (and more coarsely meshed)
  surface should be the master.
\item Problems may occur if single nodes get in contact and if
  surfaces with corners are sliding against each other. Smoothing the
  surfaces may be helpful.
\item Nodes of the master surface may penetrate the slave surface;
  again, smoothing the surfaces may reduce this, see fig.~\ref{fig:contactPenetrate}.
\begin{figure}
\includegraphics[width=0.4\textwidth]{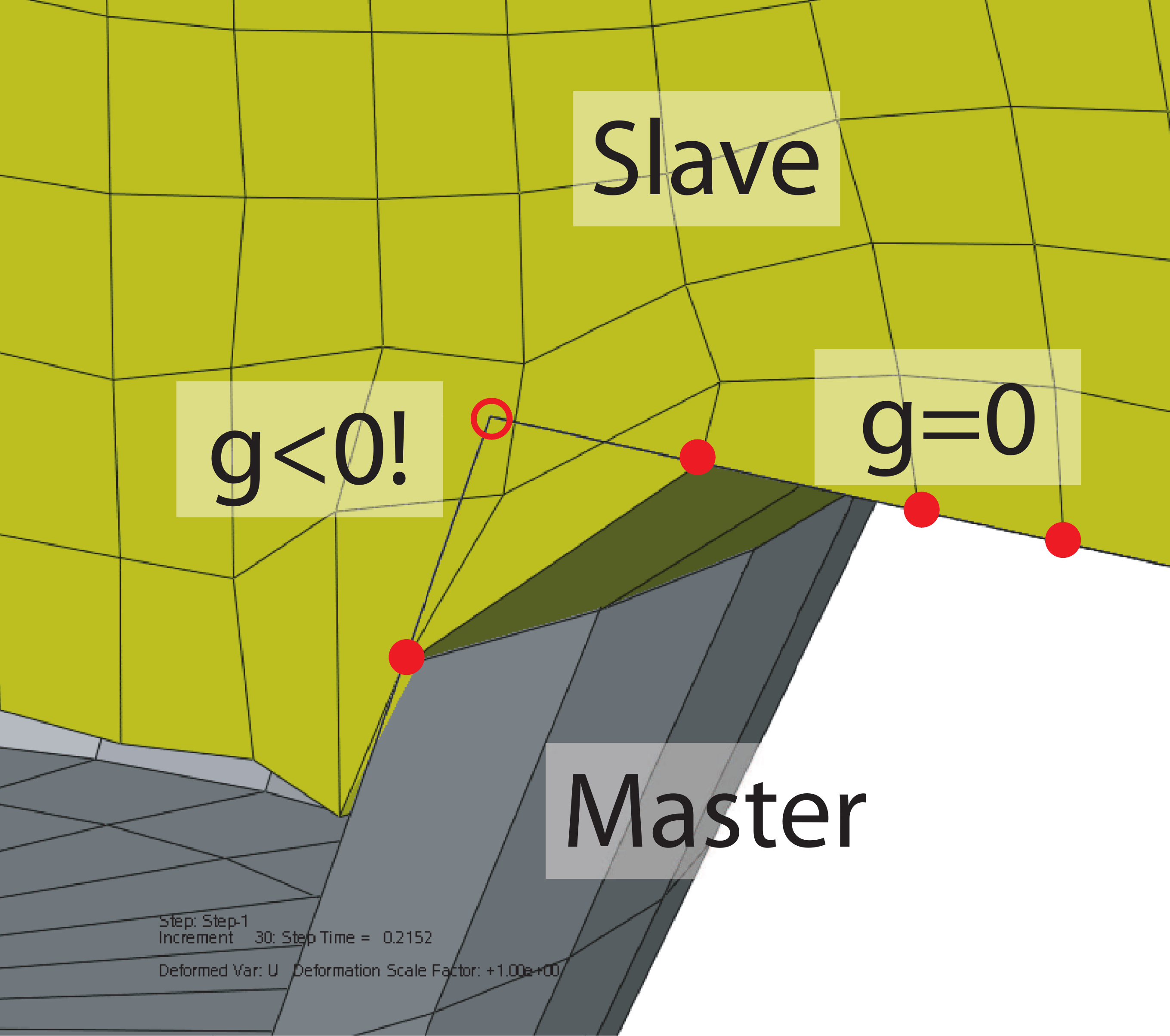}
\includegraphics[width=0.4\textwidth]{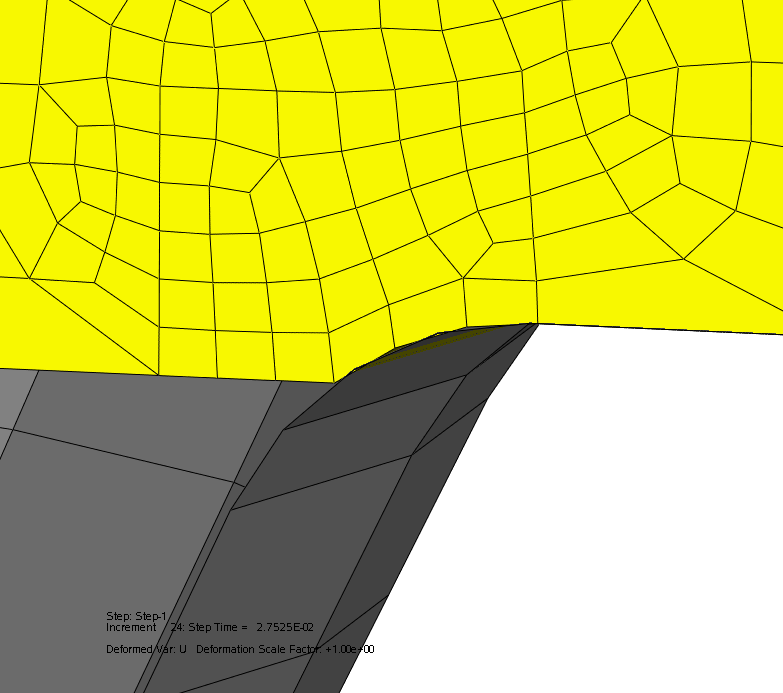}\\
\caption{Left: A node on a sharp corner on the master surface may
  penetrate the slave surface (distance of node from surface
  $g<0$). Right: Smoothing the corner and refining the slave mesh
  reduces (but does not completely eliminate) the penetration.}
\label{fig:contactPenetrate}
\end{figure}
\item Some discretization error is usually unavoidable if curved
  surfaces are in contact. With a pure master-slave algorithm,
  penetration and material overlap are the most common problem; with a
  symmetric choice (both surfaces are used as master and as slave),
  gaps may open between the surfaces, see fig.~\ref{fig:masterSlave}. Check for discretization errors
  in the postprocessor.
\begin{figure}
\begin{center}
\subfigure[Workpiece as master]{\includegraphics[width=0.32\textwidth]{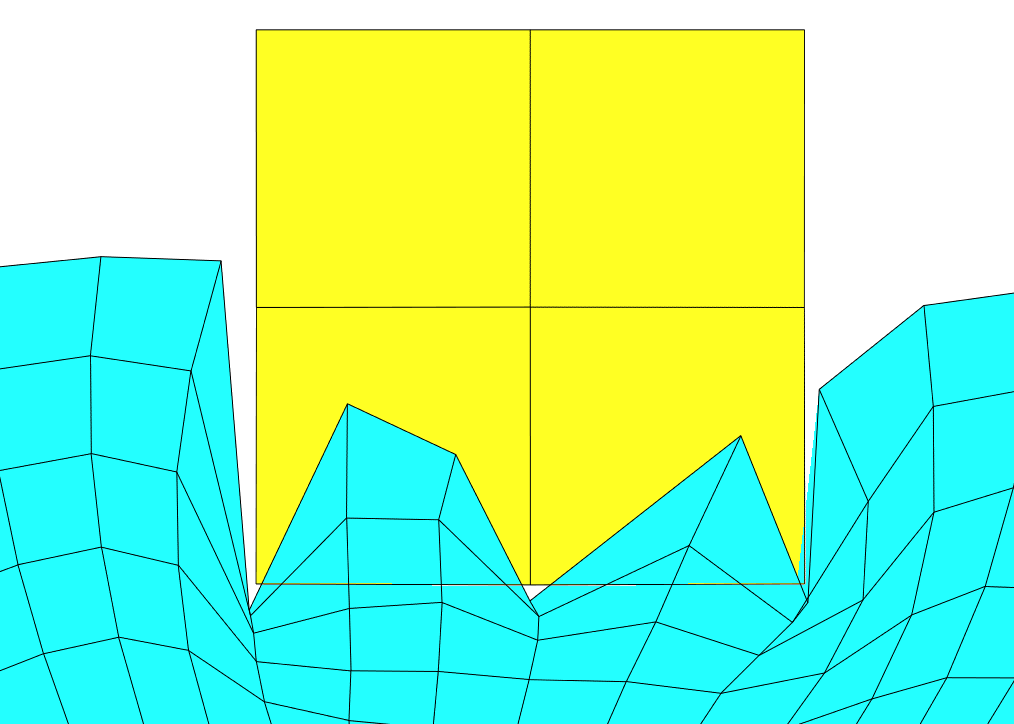}}
\subfigure[Hammer as master]{\includegraphics[width=0.32\textwidth]{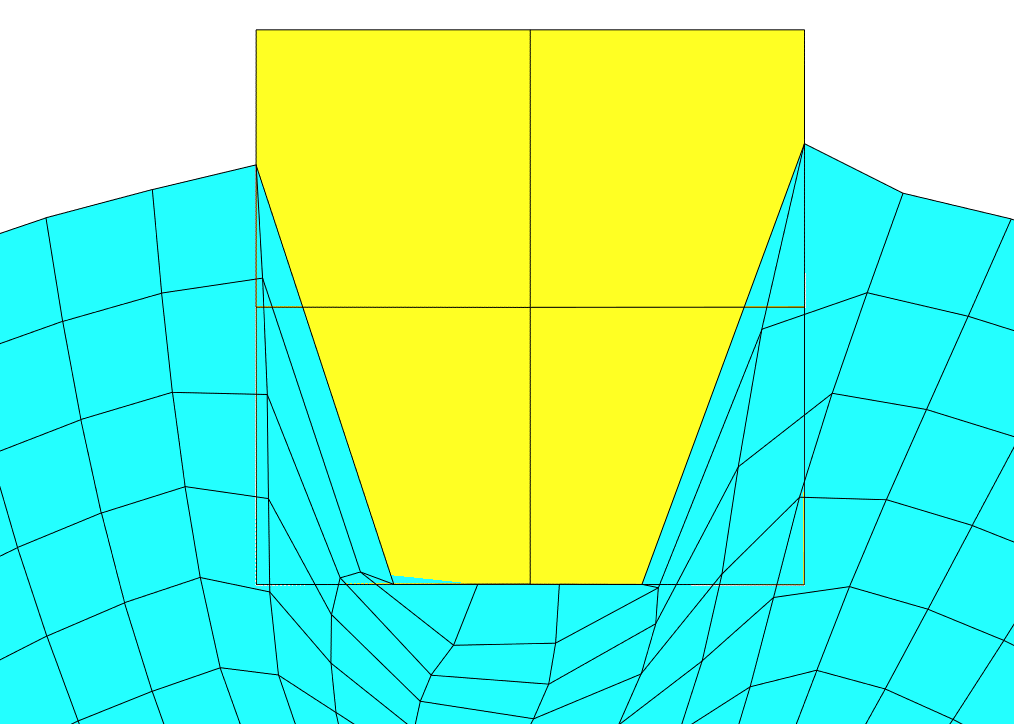}}
\subfigure[Symmetric master-slave contact]{\includegraphics[width=0.32\textwidth]{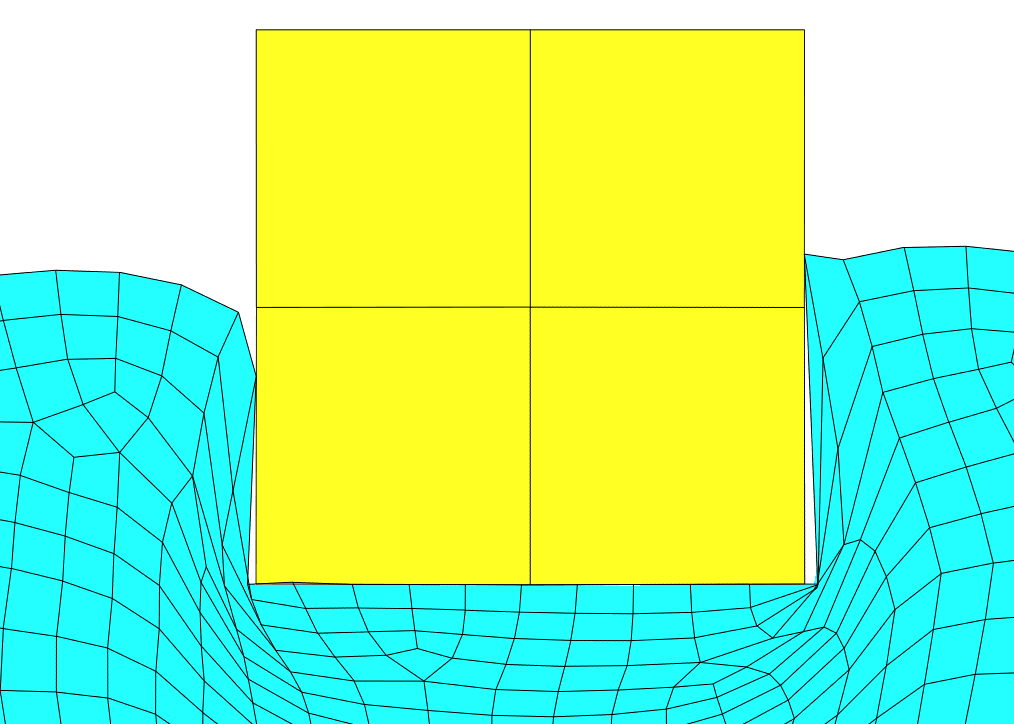}}
\end{center}
\caption{Master-slave algorithm for the example of a hammer hitting
  a workpiece. Using a master-slave algorithm results in penetrations
  of master nodes into the slave surface; a symmetric choice avoids
  this, but causes gaps to form.\label{fig:masterSlave}}
\end{figure}
\item  Discretization errors may also affect the contact
  force. Consider, for example, the Hertzian contact problem of two
  cylinders contacting each other. If the mesh is coarse, there will be a
  notable change in the contact force whenever the next node comes
  into contact. Spurious oscillations of the force may be caused by this.
\item Make sure that rigid-body motion of contact partners before the
  contact is established is removed either by adding appropriate
  constraints or by using a stabilization procedure.
\item Second-order elements may cause problems in contact (see
  \ref{item:contactSecond} and fig.~\ref{fig:nodeRelease}) \cite{Krueger2004,Baeker2018}; if they do,
  try switching to first-order elements.
\end{enumerate}

\subsection{Other considerations}
\begin{enumerate}
\item If you are inexperienced in using finite elements, start with
  simple models. Do not try to directly set up a complex model from
  scratch and make sure that you understand what your program does and
  what different options are good for. It is almost impossible to find
  errors in a large and complex model if you do not have long
  experience and if you do not know what results you expect beforehand.
\item Many parameters that are not specified by the user are set to
  default values in finite element programs. You should check whether
  these defaults are correct; especially for those parameters that
  directly affect the solution (like element types, material
  definitions etc.). If you do not know what a parameter does and
  whether the default is appropriate, consult the manual. 
For parameters that only affect the efficiency of the solution (for
example, which solution scheme is used to solve matrix equations),
understanding the parameters is less important because a wrongly
chosen parameter will not affect the final solution, but only the CPU
time or whether a solution is found at all.

\item \label{item:special}
Modern finite element software is equipped with a plethora of
  complex special techniques (XFEM, element deletion, node separation,
  adaptive error-controlled mesh-refinement, mixed Eulerian-Lagrangian
  methods, particle based methods, fluid-structure interaction,
  multi-physics, 
  user-defined subroutines etc.). If 
  you plan to use these techniques, make sure that you understand them
  and test them using simple models. If possible, build up a basic model without
these features first and then add the complex behaviour.
Keep in mind that the impressive
  simulations you see in presentations were created by experts and may have been carefully
  selected and may not be typical for the performance.
\end{enumerate}

\section{Solving the model}
Even if your model is solved without any
convergence problems, nevertheless look at the log file written by the
solver to check for warning messages. They may be harmless, but they
may indicate some problem in defining your model.

Convergence problems are usually reported by the program with warning
or error messages. 
You can also see that your model has not converged if the
final time in the time step is not the end time you specified in the
model definition.

There are two reasons for convergence problems: On the one hand, the
solution algorithm may fail to find a solution albeit a solution of
the problem does exist. In this case, modifying the solution algorithm
may solve the problem (see
section~\ref{sec:modifySolutionAlgorithm}). On the other hand, the
problem definition may be faulty so that the problem is unstable
and does not have a solution (section~\ref{sec:Instabil}).

 If you are new to finite element simulations, you
may be tempted to think that these errors are simply caused by
specifying an incorrect option or forgetting something in the model
definition. Errors of this type exist as well, but they are usually
detected before calculation of your model begins (and are not
discussed here). Instead, 
treat the non-convergence of your simulation in the same way
  as any other scientific problem. Formulate hypotheses why the
  simulation fails to converge. Modify your model to prove\footnote{Of
  course natural science is not dealing with ``proofs'', but this is
  not the place to think about the philosophy of science. Replace
  ``prove'' with ``strengthen'' or ``find evidence for'' if you like.} or disprove
these hypotheses to find the cause of the problems.

\subsection{General considerations}

\begin{enumerate}
\item\label{item:increment}
In an implicit simulation, the size of the time increments is
 usually automatically controlled by the program. If convergence is
 difficult, the time increments are reduced.\footnote{The rationale
   behind this is that the solution from the previous increment is a
   better initial guess for the next increment if the change in the
   load is reduced.} Usually, the program stops if the time increment
 is too small or if the convergence problems persist even after
 several cutbacks of the time increment. (In Abaqus, you get the error
 messages \co{Time increment smaller than minimum} or \co{Too many
   attempts}, respectively.) These messages themselves thus do not tell you
 anything about the reason for the convergence problems. 

To find the
 cause of the convergence problems, look at the solver log file
 in the increment(s) before the final error message. You will probably
 see warnings that tell you what kind of convergence problem was responsible (for
 example, the residual force is too large, the contact algorithm did
 not converge, the temperature increments were too large).
If available, also look at the unconverged solution and compare it to
the last, converged timestep. Frequently, large changes in some
quantity may indicate the location of the problem.
\item Use the postprocessor to identify the node with the largest
  residual force and the largest change in displacement in the final
  increment. Often (but not always) this tells you where the problem
  in the model occurs. (Apply the same logic in a thermal simulation
  looking at the temperature changes and heat fluxes.)
\item If the first increment does not converge, set the size of the 
  first time increment to a very small value. If the problem persist,
  the model itself may be unstable (missing boundary conditions,
  initial overlap of contacting surfaces). To find the cause of the
  problem, you can remove all external loads step by step or add
  further boundary conditions to make sure that the model is properly
  constrained (if you pin two nodes for each component, rigid body
  movements should be suppressed~-- if the model converges in this
  case, you probably did not have sufficient boundary conditions in
  your original model). Alternatively or additionally, you may add numerical
  stabilization to the problem definition. (In numerical
  stabilization, artificial friction is added to the movement of nodes
  so that stabilizing forces are generated if nodes start to move
  rapidly.) However, make sure that the stabilization does not affect
  your results too strongly.

  Also check for abrupt jumps
  in some boundary conditions, for example a finite displacement that
  is defined at the beginning of a step or a sudden jump in
  temperature or load. If you apply a load instantaneously, cutting
  back the time increments does not help the solution process. If this
  occurs, ramp your load instead.
\item\label{item:rapidChange} Avoid rapid changes in an amplitude within
  a calculation step (see also \ref{item:timeSteps} and \ref{item:timestep2}). For
  example, if you hold a heat flux (or temperature or stress) for a
  long time and then abruptly 
  reduce it within the same calculation step, the time increment will
  suddenly jump to a point where the temperature is strongly
  reduced. This abrupt change may cause convergence problems. Define a
  second step and 
choose small increments at the beginning of the
  second step where large changes in the model can be expected.
\item Try the methods described in
  section~\ref{sec:modifySolutionAlgorithm} to see whether the problem
  can be resolved by changing the solution algorithm.
\item Sometimes, it is the calculation of the material law at an
  integration point that does not converge (to calculate stresses from
  strains at integration point inside the solver, another Newton
  algorithm is used at each integration point
  \cite{baker2002numerische}). If this is the case, the material
  definition may be incorrect or problematic (for example, due to incorrectly
  specified material parameters or because there is extreme softening
  at a point).
\item\label{item:minimum} Simplify your model step by step to find the reason of the
  convergence problems. Use simpler material laws (simple plasticity
  instead of damage, elasticity instead of plasticity), switch off
  non-linear geometry, remove external loads etc. 
 If the problem persists, try to create a minimum example~-- the
 smallest example you can find that shows the same problem. This has
 several advantages: the minimum example is easier to analyse, needs
 less computing time so that trying things is faster, and it can also
 be shown to others if you are looking for help (see section~\ref{sec:help}).
\item If your simulation is static, switching to an implicit dynamic simulation
  may help because the inertial forces act as natural stabilizers. If
  possible, use a quasi-static option. 
\item \label{it:explicitConvergence}
Explicit simulations usually have less convergence problems. A
  frequently-heard advice to solve convergence problems is to switch
  from implicit to explicit models. 
I strongly recommend to only switch from
  implicit static to explicit quasi-static for convergence reasons if
  you understand the reasons of the convergence problems and cannot
  overcome them with the techniques described here. You should also keep in mind
  that explicit programs may offer a different functionality (for
  example, different element types). If your problem is static, you
  can only use a quasi-static explicit analysis which may also have
  problems (see~\ref{it:explicitScaling}). Be aware that in an
  explicit simulations, elastic waves may occur that may change the
  stress patterns.
\end{enumerate}

\subsection{Modifying the solution algorithm}
\label{sec:modifySolutionAlgorithm}

If your solution algorithm does not converge for numerical reasons,
these modifications may help. They are useless if there is a true
model instability, see section~\ref{sec:Instabil}.

\begin{enumerate}
\item Finite element programs use default values to control the Newton
  iterations. If no convergence is reached after a fixed number of
  iterations, the time step is cut back. In strongly non-linear
  problems, these default values may be too tight. For example, Abaqus
  cuts back on the time increment if the Newton algorithm does not
  converge after 4~iterations; setting this number to a larger value
is often sufficient to reach convergence (for example, by adding \co{*Controls,
  analysis=discontinuous} to the input file).
\item If the Newton algorithm does not converge, the time increment is
  cut back. If it becomes smaller than a pre-defined minimum value,
  the simulation stopswith an error message.  This minimum size of the time increment can
  be adjusted. Furthermore, if a sudden loss in stability (or
  change in load)  occurs so
  that time increments need to be changed by several orders of
  magnitude, the number of cutbacks also needs to be adapted (see next
  point). In this case, another option is to define a new time step (see \ref{item:rapidChange})
  that starts at this critical point and that has a small initial increment.
\item
The allowed number of cutbacks 
  increment can also be adapted (in Abaqus, use \co{*CONTROLS,
    parameters=time incrementation}). This may be helpful if the
  simulation proceeds at first with large increments before some
  difficulty is reached -- allowing for a larger number of cutbacks
  enables the program to use large timesteps at the
  beginning. Alternatively, you can reduce the maximum time increment
  (so that the size of the necessary cutback is reduced)
  or you can split your simulation step in two with different time
  incrementation settings in the step where the problem occurs  (see
  \ref{item:rapidChange}).  
\item Be aware that the previous two points will work sometimes, but
  not always. There
  is usually no sense in allowing a smallest time increment that is
  ten or twenty orders of magnitude smaller than the step size or to
  allow for dozens of cutbacks, this only increases the CPU time.
\item Depending on your finite element software, there may be many
  more options to tune the solution process. In Abaqus, for example,
  the initial guess for the solution of a time increment is calculated
  by extrapolation from the previous steps. Usually this improves
  convergence, but it may cause problems if something in the model
  changes abruptly. In this case, you can switch the extrapolation
  off (\co{*STEP, extrapolation=no}). You can also add a line search algorithm that scales the
  calculated displacements to find a better solution (\co{*CONTROLS,
    parameters=line search}). Consult the manual for options to
  improve convergence.
\item While changing the iteration control (as explained in the
  previous points) is often needed to achieve convergence, 
  the solution controls that are used to determine whether a solution
  has converged should only be changed if absolutely necessary. Only
  do so 
(in Abaqus, use \co{*CONTROLS,
    parameters=field})
if you know exactly what you are doing. 
One example where changing the controls may be necessary is when the stress is strongly concentrated in a small part of a very large structure \cite{baker2002}. In this
case, an average nodal force that is used to determine convergence may
impose too strong a constraint on the convergence of the solution, so
that convergence should be based on local forces in the region of
stress concentration. Be aware that since forces, not stresses, are
used in determining the convergence, changing the mesh density
requires changing the solution controls.
  
Make
  sure that the accepted solution is indeed a solution and that your
  controls are sufficiently strict. Vary the controls to ensure that
  their value does not affect the solution.
\item Contact problems sometimes do not converge due to problems in
  establishing which nodes are in contact (sometimes called
  ``zig-zagging'' \cite{wriggers2006}). This often happens if the
  first contact is made by
  a single node. Smoothing the contact surfaces may help.
\item If available and possible, use general contact definitions where
  the contact surfaces are determined automatically.
\item \label{it:softContact}
If standard contact algorithms do not converge, soft contact
  formulations (which implement a soft transition between ``no
  contact'' and ``full contact'') may improve convergence; however,
  they may allow for some penetration of the surfaces and thus affect
  the results.

\end{enumerate}

\subsection{Finding model instabilities}
\label{sec:Instabil}
A model is unstable if there actually is no solution to the mechanical
problem. 
\begin{enumerate}
\item Instabilities are frequently due to a loss in load bearing
  capacity of the structure. There are several reasons for that: 
\begin{itemize}
\item The material definition may be incorrect. If, for example,  a plastic material is defined
  without hardening, the load cannot increase after the component has
  fully plastified. Simple typos or incorrectly used units may also
  cause a loss in material strength. 
\item Thermal softening (the reduction of strength with increasing
  temperature) may cause an instability in a thermo-mechanical problem.
\item  Non-linear geometry may cause an instability because the cross
  section of a load-bearing component reduces during deformation.
\item A change in contact area, a change from sticking to sliding in a
  simulation with friction or a complete loss of contact between
  to bodies may also cause instabilities because the structure may not
  be able to bear an increase in the load.
\end{itemize}
\item Local instabilities may cause highly distorted
  meshes that prevent convergence. It may be helpful to define the
  mesh in such a way that elements become more regular during
  deformation (see also~\ref{it:distort}). 
\item If your model is load-controlled (a force is applied), switch to
  a displacement-controlled loading. This avoids instabilities due to
  loss in load-bearing capacity. 
\item Artificial damping (stabilization) may be added to stabilize an unstable
  model. However, check carefully that the solution is not unduly
  affected by this. Adding artificial damping may also help to
  determine the cause of the instability. If your model converges with
  damping, you know that an instability is present.
\end{enumerate}

\subsection{Problems in explicit simulations}
As already stated in \ref{it:explicitConvergence}, explicit
simulations have less convergence problems than implicit
simulations. 
However, sometimes even an explicit simulation may run into trouble.

\begin{enumerate}
\item During simulation, elements may distort excessively. This may
  happen for example if a concentrated load acts on a node or if the
  displacment of a node becomes very large due to a loss in
  stability (for example in a damge model). In this case, the element
  shape might become invalid (crossing over of element edges, negative
  volumes at integration points etc.). If this happens, changing the
  mesh might help -- elements that have a low quality (large aspect
  ration, small initial volume) are especially prone to this type of
  problem. Note that second-order elements are often more sensitive to
  this problem than first-order elements.
\item The stable time increment in an explicit simulation is given by
  the time a sound wave needs to travel through the smallest
  element. If elements distort strongly, they may become very thin in
  one direction so that the stable time increment becomes unreasonably
  small. In this case, changing the mesh might help.
\end{enumerate}

\section{Postprocessing}
There are two aspects to checking that a model is correct: Verification
is the process of showing that the model was correctly specified and
actually does what it was created to do (loads, boundary conditions,
material behaviour etc.\ are correct). Validation means to check the
model by making an independent prediction (i.\,e., a prediction that
was not used in specifying or calibrating the model) and checking this prediction
in some other way (for example, experimentally).\footnote{Note that the terms
  ``verification'' and ``validation'' are used differently in
  different fields.}

\textbf{General advice}
If you modify your model significantly (because you build up a
complicated model in steps, have to correct
errors or add more complex material behaviour to get agreement with
experimental results etc.), you should again check the model. It is
not clear that the mesh density that was sufficient for your initial
model is still sufficient for the modified model. The same is true for
other considerations (like the choice of element type etc.)

\subsection{Checking the plausibility and verifying the model}\label{sec:validierung}
\begin{enumerate}
\item Check the plausibility of your results. If your simulation deviates from your intuition,
  continue checking until you are sure that you understand why your
  intuition (or the simulation) was incorrect. \emph{Never} believe a
  result of a simulation that you do not understand and that should be
  different according to your intuition. Either the model or your
  understanding of the physical problem is incorrect -- in both cases,
  it is important to understand all effects.
\item Check your explanations for the solution, possibly with
  additional simulations. For example, if you assume that thermal
  expansion is
  the cause of a local stress maximum, re-run the simulation with a
  different or vanishing coefficient of thermal expansion. Predict the
  results of such a simulation and check whether your prediction was
  correct. 
\item Check all important solution variables. Even if you are only
  interested in, for example, displacements of a certain point, check
   stresses and strains throughout the model.
\item In 3D-simulations, do not only look at contour plots of the
  component's surface; also check the results inside the component by
  cutting through it.
\item Make sure you understand which properties are vectors or
  tensors. Which component of stresses or strains are relevant depends
  on your model, the material, and the question you are trying to
  answer. Default settings of the postprocessor are not always
  appropriate, for example, Abaqus plots the von-Mises-stress as
  default stress variable, which is not very helpful for ceramic materials.
\item Check the boundary conditions again. Are all nodes constrained
  in the desired manner?
Exaggerating the deformation (use \co{Common plot options} in Abaqus)
or picking nodes with the mouse
may be helpful to check this precisely. 
\item Check the mesh density (see~\ref{item:meshconverge}). If possible, calculate the model with
  different mesh densities (possibly for a simplified problem) and
  make sure that the mesh you finally use is sufficiently fine. When
  comparing different meshes, the variation in the mesh density should
  be sufficiently large to make sure that you can actually see an
  effect.
\item Check the mesh quality again, paying special attention on
  regions where gradients are large. Check that the conditions
  explained in section~\ref{sec:mesh} (element shapes and sizes, no strong
   discontinuities in the element sizes) are fulfilled and that
   discontinuities in the stresses are not due to a change in the
   numerical stiffness (due to a change in the integration scheme or
   element size).
\item Check that stresses are continuous between elements. At
  interfaces between different materials, check that normal stresses
  and tangential strains are continuous.
\item Check that the normal stress at any free surface is zero.
\item Check the mesh density at contact surfaces: can the actual
  movement and deformation of the surfaces be represented by the mesh?
  For example, if a mesh is too coarse, nodes may be captured in a
  corner or a surface may not be able to deform correctly. 
\item  Keep in mind that discretization errors at contact surfaces
  also influence stresses and strains. If you use non-standard contact
  definitions (\ref{it:softContact}), try to evaluate how these
  influence the stresses (for example by comparing actual node
  positions with what you would expect for hard contact).
\item Watch out for divergencies. The stress at a sharp notch or cack
  tip is theoretically infinite~-- the value shown by your program is
  then solely determined by the mesh density and, if you use a contour
  plot, by the extrapolation used by the postprocessor (see~\ref{item:extrapolate}).
\item In dynamic simulations, elastic waves propagate through the
  structure. They may dominate the stress field. Watch out for
  reflections of elastic waves and keep in mind that, in reality,
  these waves are dampened.
\item If you assumed linear geometry, check whether strains and
  deformations are
  sufficiently small to justify this assumption, see~\ref{item:nlgeom}.
\smallskip

\end{enumerate}

\subsection{Implementation issues}
\begin{enumerate}
\item\label{item:extrapolate} Quantities like stresses or strains are only defined at
  integration points. Do not rely on extreme values from a contour
  plot~-- these values are extrapolated. It strongly depends on the
  problem whether these extrapolated values are accurate or not. For
  example, in an elastic material, the extrapolation is usually
  reasonable, in an
ideally-plastic material, extrapolated von Mises
  stresses may exceed the actual yield stress by a factor
  of~2 or more. Furthermore, the contour lines themselves may show incorrect maxima
  or minima, see
  fig.~\ref{fig:extrapolate} for an example.
\begin{figure}
\includegraphics[width=0.48\textwidth]{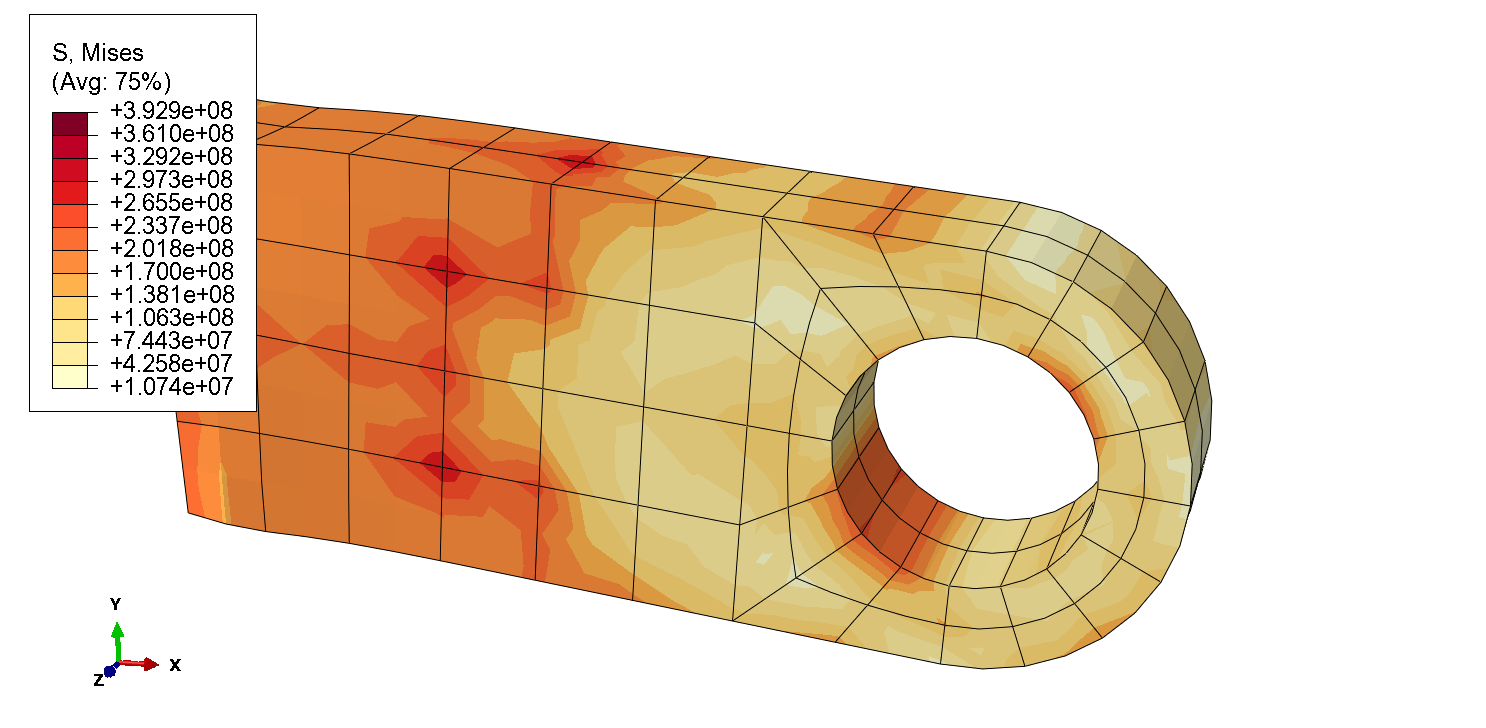}
\quad
\includegraphics[width=0.48\textwidth]{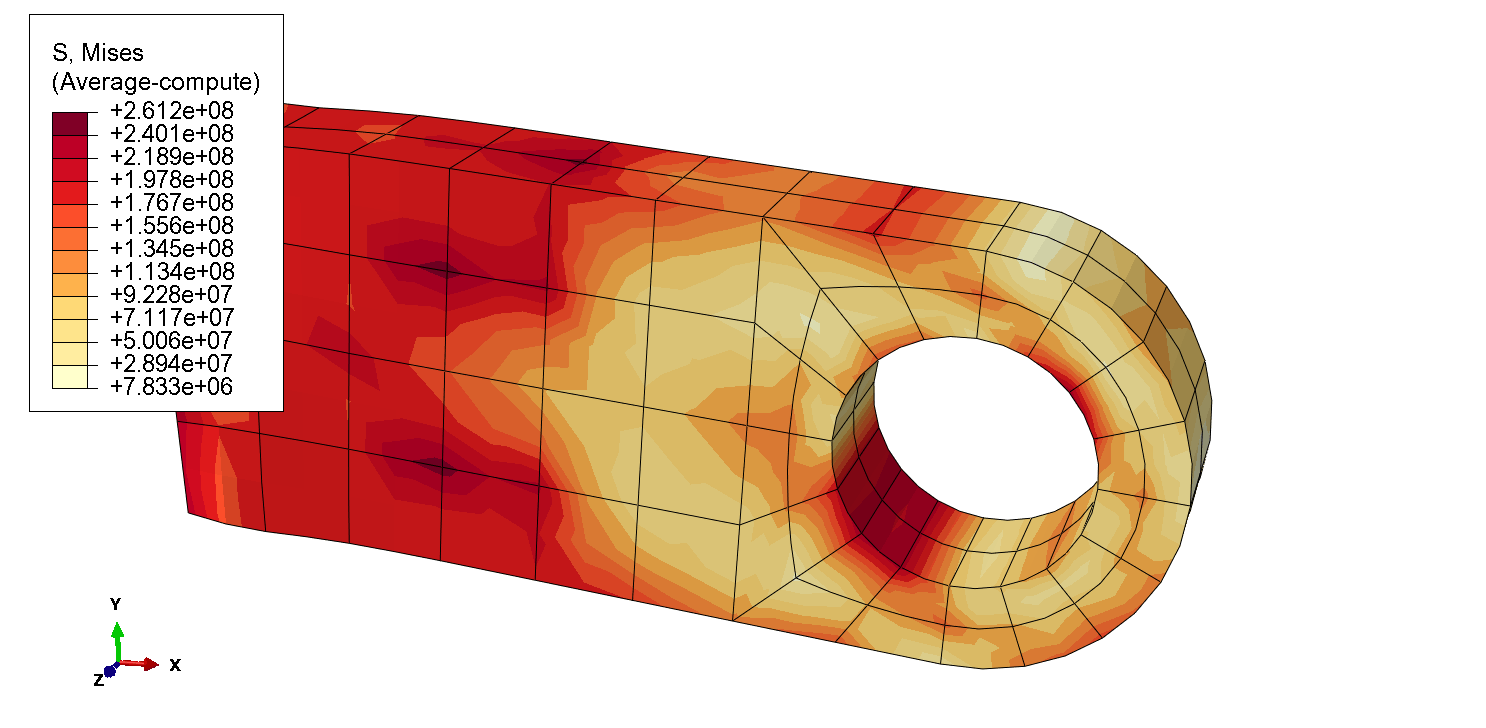}

\includegraphics[width=0.48\textwidth]{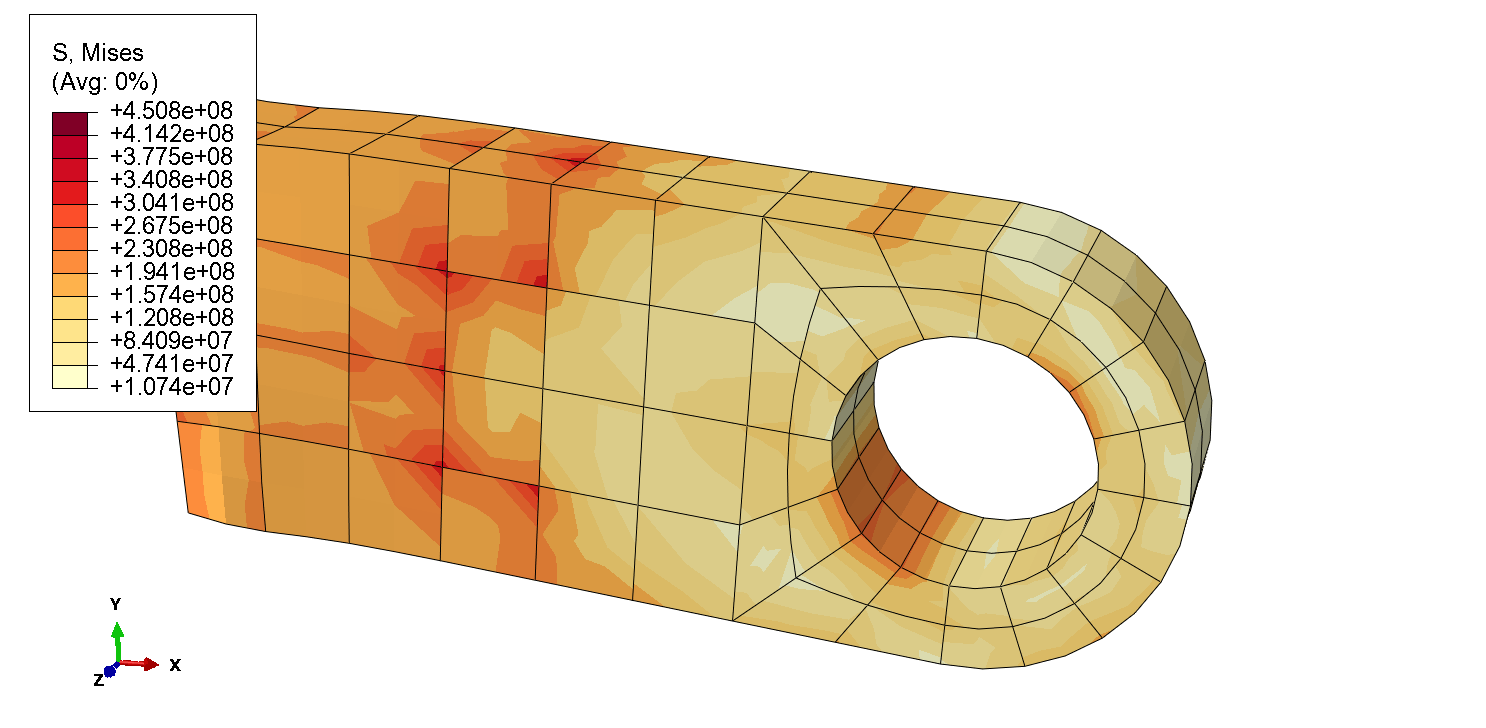}
\quad
\includegraphics[width=0.48\textwidth]{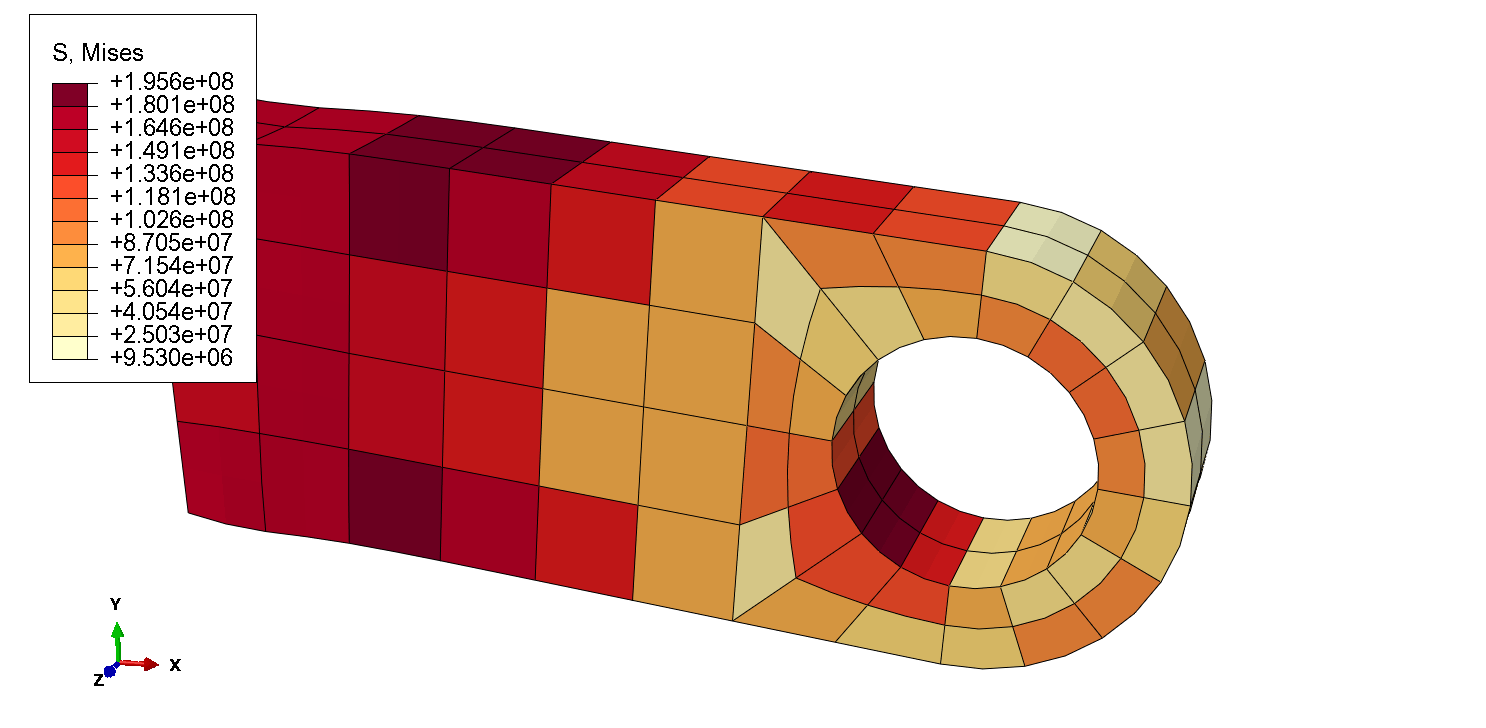}
\caption{Von Mises stress in a simple lug constrained on the left and loaded in the hole
  on the right. The material is ideally plastic with a yield stress of
\numprint[MPa]{180} so that the von Mises stress at the integration
points never exceeds this value. All four figures show the same
simulation result. Due to extrapolation from the integration points,
the maximum value in the contour plot is much too large, except for
the quilt plot on the lower right. The exact
maximum value depends on  how values are
extrapolated and averaged at the surface and between elements. }
\label{fig:extrapolate}
\end{figure}
\item It is often helpful to use ``quilt'' plots where each element is
  shown in a single color averaged from the integration point values
  (see also fig.~\ref{fig:extrapolate}).
\item The frequently used rainbow color spectrum has been shown to be
  misleading and should not be used \cite{borland2007}. Gradients may be difficult to
  interpret because human color vision has a different sensitivity in
  different parts of the spectrum. Furthermore, many people have a
  color vision deficiency and are unable to discern reds, greens and
  yellows. For variables that run from zero to a maximum value
  (temperature, von-Mises stress), use a
  sequential spectrum (for example, from black to red to yellow), for
  variables that can be positive and negative, use a diverging
  spectrum with a neutral color at zero, see fig.~\ref{fig:rainbow}.
\begin{figure}
\includegraphics[width=\textwidth]{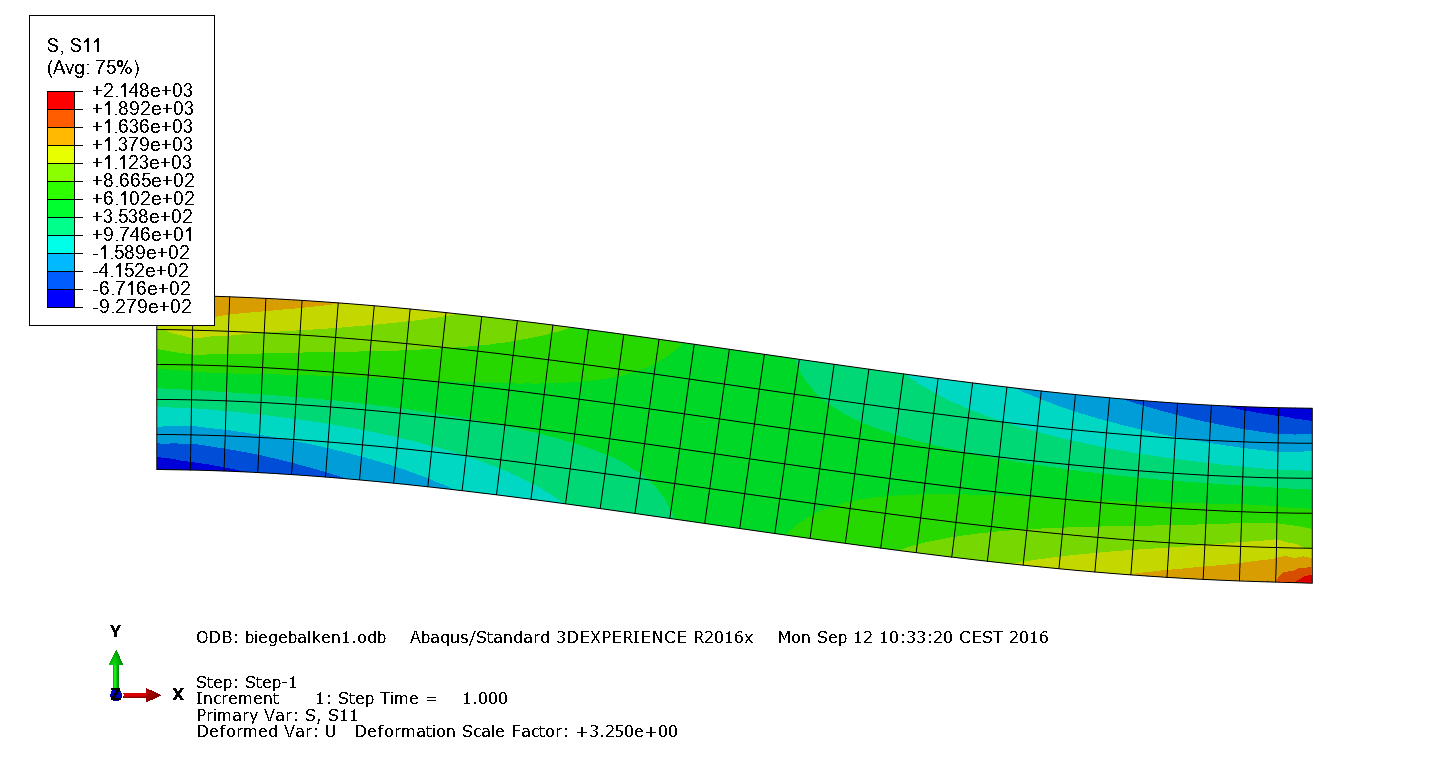}

\includegraphics[width=\textwidth]{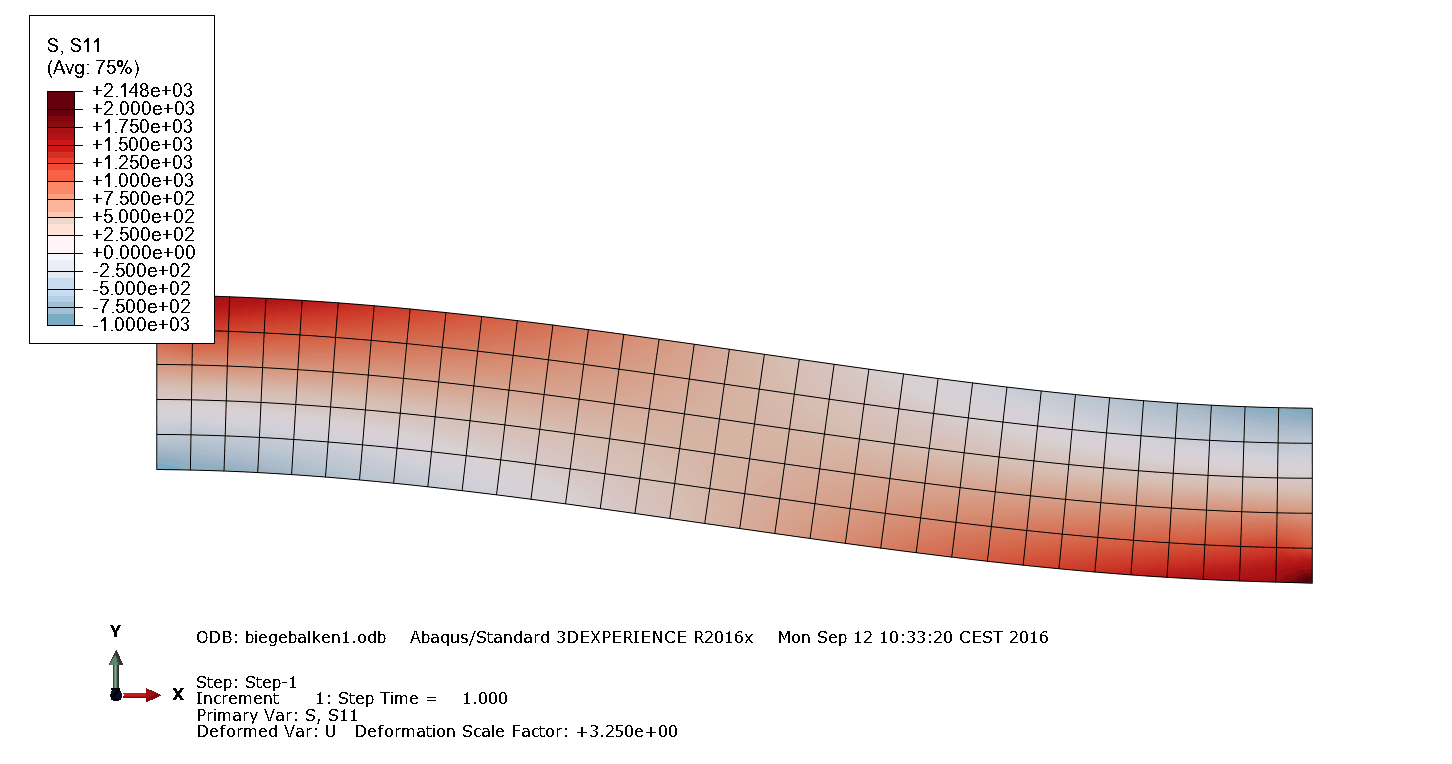}
\caption{Plot of the stress component $\sigma_{11}$ in a bar under
  tension and bending. In a rainbow plot, it is difficult to see
  whether the gradient is homogeneous or not and to find the line of
  zero stress. A divergent spectrum with a neutral color at zero makes
  understanding the stress field easier. Furthermore, it is easier to
  interpret for people with color vision deficiency. The plot was done
using the Abaqus plugin SpectrumBaker.}
\label{fig:rainbow} 
\end{figure}
\item Discrete time-stepping (see~\ref{item:timeSteps}) may also
  influence the post-processing of results. If 
  you plot the stress-strain curve of a material point by connecting
  values measured at the discrete simulation times, the resulting
  curve will not coincide perfectly with the true stress-strain although
  the data points themselves are correct.
\item Complex simulation techniques (like XFEM, element deletion etc.,
  see~\ref{item:special}) frequently use internal parameters to
  control the simulation that may affect the solution process. Do not
  rely on default values for these parameters and check that the
  values do not affect the solution inappropriately.
\item If you use element deletion, be aware that removing elements from the simulation is
  basically an unphysical process since material is removed. This may
  affect the energy balance or stress fields near the removed
  elements.  For example, in models of machining
  processes, removing elements at the tool tip to separate the
  material strongly influences the residual stress field.
\end{enumerate}

\subsection{Validation}\label{sec:verifikation}

\begin{enumerate}
\item If possible, use your model to make an independent prediction
  that can be tested. 
\item If you used experimental data to adapt unknown parameters
  (see~\ref{sec:eingangsgroessen}), correctly reproducing these data 
  with the model does not validate it, but only verifies it.
\item The previous point also holds if you made a prediction and
  afterwards had to change your model to get agreement with an
  experiment. After this model change, the experiment cannot be
  considered an independent verification.
\end{enumerate}

\section{Getting help}\label{sec:help}
If you cannot solve your problem, you can try to get help from the
support of your software (provided you are entitled to support) or
also from the internet (for 
example on researchgate or imechanica). To get helpful answers, please observe the following points:
\renewcommand*{\theenumi}{\thesection-\arabic{enumi}}
\begin{enumerate} 
\item Check that you have read relevant pages in the manual and that
  your question is not answered there.
\item Describe your problem as precisely as possible. Which error did
  occur? What was the exact error message and which warnings did
  occur? Show pictures of the model and describe 
  the model (which element type, which material, what kind of
  problem~-- static, dynamic, explicit, implicit etc.).
\item If possible, provide a copy of your model or, even better,
  provide a minimum example that shows the problem (see~\ref{item:minimum}).

\item If you get answers to your request, give feedback
  whether this has solved your problem, especially if you are in an
  internet forum or similar. People are
  sacrificing their time to help you and will be interested to see
  whether their advice was actually helpful and what the solution to
  the problem was. Providing feedback will also help others who find
  your post because they are facing similar problems.
\end{enumerate}

\section*{Acknowledgement}
Thanks to Philipp Seiler for many discussions and for reading a draft
version of this manuscript, and to Axel Reichert for sharing his
experience on getting models to converge. 

\bibliography{HowToLiterature}

\end{document}